\documentclass[reqno,11pt]{article}
\usepackage{amsmath,amssymb}
\usepackage{amsthm,amscd,amsfonts,}
\usepackage{amsmath}
\usepackage{amsthm}
\usepackage{graphicx}
\usepackage{amsmath,amsthm}
\usepackage{amsmath,amstext,amsfonts}
\usepackage{subfigure}
\usepackage[T1]{fontenc}
\usepackage[utf8]{inputenc}
\usepackage{authblk}

\newtheorem{theorem}{Theorem}[section]

\newtheorem{corollary}{Corollary}[section]
\newtheorem{example}{Example}[section]

\newcommand{\R}{\mathbb{R}}
\numberwithin{equation}{section}

\begin{document}

\title{\textbf{Integrating Factor and First Integrals of a Class of Third Order Differential Equations}}
\author{M. Al-Jararha\thanks{mohammad.ja@yu.edu.jo}}
\affil{Department of Mathematics, Yarmouk University, Irbid, Jordan, 21163.}

\date{}
\maketitle

\vspace{0.2in}

\noindent \textbf{AMS (2000) Subject Classification}: 34A25, 34A30.

\begin{abstract}
The principle of finding an integrating factor for a none exact differential equations is extended to a class of third order differential equations. If the third order equation is not exact, under certain conditions,  an integrating factor exists which  transforms it to an exact one. Hence, it can be reduced into a second order differential equation. In this paper, we give explicit forms for certain integrating factors of a class of the third order differential equations.
\end{abstract}
\maketitle
\section{Introduction}

Third order nonlinear differential equations play a major role in Applied Mathematics, Physics, and Engineering \cite{Ames, Jordan, Lefschetz, Struble,von}. To find the general solution of a third order nonlinear differential equation is not an easy problem in the general case. In fact, a very specific class of nonlinear third order differential equations can be solved by using special transformations. Other  technique is to reduce the order of the differential equation into the second order, by finding a proper integrating factor. Recently, many studies appear to deal with the problem of the existence of an integrating factor for certain differential equations. For example, in \cite{AlAhmad, existencepp, Cheb, Bouquet}, the authors investigated the existence of an integrating factor of some classes of second order differential equations. In \cite{existencepp}, the authors investigated the existence of an integrating factor of $n-th$ order differential equations which has known symmetries of certain type. 

In \cite{kai}, the authors improve some symbolic algorithms to compute the integrating factor for certain class of third order nonlinear differential equation. In this paper, we investigate the existence of an integrating factor of the following class of third order nonlinear differential equations: 
\begin{equation}\label{thirdorder1}
F_3(t,y,y^\prime,y^{\prime\prime})y^{\prime\prime\prime}+F_2(t,y,y^\prime,y^{\prime\prime})y^{\prime\prime}+F_1(t,y,y^\prime,y^{\prime\prime})y^\prime+F_0(t,y,y^\prime,y^{\prime\prime})=0.
\end{equation}
In fact, we present some theoretical results related to the existence of certain forms of the integrating factor for \eqref{thirdorder1}. We also present some illustrative examples.  
\section{Integrating Factor and First Integral of a Class of Third Order Differential Equations}
In this section, we investigate the existence of special forms of integrating factor of a class of none exact third order differential equations. In general, we said that the $n-th$ order differential equation $f\left(t,y,y^\prime,\cdots,y^{(n-1)},y^{(n)}\right)=0$ is exact if there exists a differentiable function $\Psi\left(t,y,y^\prime,\cdots,y^{(n)}\right)=c$, such that $\frac{d}{dt}\Psi\left(t,y,y^\prime,\cdots,y^{(n)}\right)=f\left(t,y,y^\prime,\cdots,y^{(n-1)},y^{(n)}\right)=0$. In this case,  $\Psi\left(t,y,y^\prime,\cdots,y^{(n)}\right)=c$ is called the first integral of the differential equation $f\left(t,y,y^\prime,\cdots,y^{(n-1)},y^{(n)}\right)=0$  (see for example \cite{cov,Peter}). In \cite{Aljararha}, the author gave the conditions so that the first integral of 
\begin{align}\label{HODE}
&F_n\left(t,y,y^\prime,y^{\prime\prime},\ldots,y^{(n-1)}\right)y^{(n)}+F_{n-1}\left(t,y,y^\prime,y^{\prime\prime},\ldots,y^{(n-1)}\right)y^{(n-1)}+\cdots \nonumber\\
&+F_{1}\left(t,y,y^\prime,y^{\prime\prime},\ldots,y^{(n-1)}\right)y^{\prime}+F_{0}\left(t,y,y^\prime,y^{\prime\prime}\ldots,y^{(n-1)}\right)=0.
\end{align}
exists. Also, he gave an explicit formula for $\Psi\left(t,y,y^\prime,\cdots,y^{(n)}\right)$. 
for particular case, we consider the following class of third order differential equation:
\begin{equation}\label{thirdorder}
F_3(t,y,y^\prime,y^{\prime\prime})y^{\prime\prime\prime}+F_2(t,y,y^\prime,y^{\prime\prime})y^{\prime\prime}+F_1(t,y,y^\prime,y^{\prime\prime})y^\prime+F_0(t,y,y^\prime,y^{\prime\prime})=0,
\end{equation}
where $F_0,F_1,F_2$ and $F_3$ are continuous with their first partial derivatives with respect to $t,y,y^\prime,$ and $y^{\prime\prime}$ on some simply connected domain $\mathcal R$ in $\R^4$. According to \cite{Aljararha}, this third order differential equation  is exact if the following conditions:
\begin{equation}\label{exactcond}
\partial_{y^{\prime\prime}}F_0=\partial_tF_3,\;\partial_{y^{\prime\prime}}
F_1=\partial_yF_3,\; \partial_{y^{\prime\prime}}F_2=\partial_{y^\prime}F_3,
\partial_{y^\prime}F_0=\partial_tF_2,\;\partial_{y^\prime}F_1=\partial_yF_2,\; 
\text{and}\;\partial_yF_0=\partial_tF_1
\end{equation}
hold. Moreover, its first integral is given by
\begin{eqnarray}\label{psithirdorder}
 \Psi\left(t,y,y^\prime,y^{\prime\prime}\right)&=& \int_{t_0}^{t}F_0\left(\xi,y,y^\prime,y^{\prime\prime}\right)d\xi+  \int_{y_0}^{y}F_1\left(t_0,\xi,y^\prime,y^{\prime\prime}\right)d\xi\nonumber\\
&+&\int_{y^\prime_0}^{y^\prime}F_2\left(t_0,y_0,\xi,y^{\prime\prime}\right)d\xi +\int_{y_0^{\prime\prime}}^{y^{\prime\prime}}F_3\left(t_0,y_0,y^\prime_0,\xi\right)d\xi \nonumber\\
&=&c,
\end{eqnarray}
for some constant $c$. Assume that \eqref{thirdorder} is not exact differential equation, then an integrating factor $\mu(t,y,y^\prime,y^{\prime\prime})$ if exists will transform equation it into an exact differential. According to the conditions \eqref{exactcond}, the integrating factor $\mu(t,y,y^\prime,y^{\prime\prime})$ of \eqref{thirdorder} must solve the following system of partial differential equations:
\begin{equation}\label{sys1}
\left\{
\begin{array}{ll}
\mu(\mathbf y)F_{3t}(\mathbf y)+\mu_t(\mathbf y)F_{3}(\mathbf y )=\mu(\mathbf y )F_{0y^{\prime\prime}}(\mathbf y)+\mu_{y^{\prime\prime}}(\mathbf y )F_{0}(\mathbf y),\\
\mu(\mathbf y)F_{2t}(\mathbf y)+\mu_t(\mathbf y)F_{2}(\mathbf y )=\mu(\mathbf y )F_{0y^\prime}(\mathbf y)+\mu_{y^\prime}(\mathbf y )F_{0}(\mathbf y),\\
\mu(\mathbf y)F_{1t}(\mathbf y)+\mu_t(\mathbf y)F_{1}(\mathbf y )=\mu(\mathbf y )F_{0y}(\mathbf y)+\mu_y(\mathbf y )F_{0}(\mathbf y),\\
\mu(\mathbf y)F_{1y^\prime}(\mathbf y)+\mu_{y^\prime}(\mathbf y)F_{1}(\mathbf y )=\mu(\mathbf y )F_{2y}(\mathbf y)+\mu_y(\mathbf y )F_{2}(\mathbf y),\\
\mu(\mathbf y)F_{1y^{\prime\prime}}(\mathbf y)+\mu_{y^{\prime\prime}}(\mathbf y)F_{1}(\mathbf y )=\mu(\mathbf y )F_{2y}(\mathbf y)+\mu_y(\mathbf y )F_{2}(\mathbf y),\\
\mu(\mathbf y)F_{2y^{\prime\prime}}(\mathbf y)+\mu_{y^{\prime\prime}}(\mathbf y)F_{2}(\mathbf y )=\mu(\mathbf y )F_{3y^\prime}(\mathbf y)+\mu_{y^\prime}(\mathbf y )F_{3}(\mathbf y),
\end{array}
\right.
\end{equation}
where $\mathbf y=(t,y,y^\prime,y^{\prime\prime})$.
To solve such system of partial differential equations is not easy. Therefore, to find an integrating factor of \eqref{thirdorder}, we look for integrating factors of certain forms. Particularly, we are looking for an integrating factor of the form $\mu(\xi)$, where $\xi:=\xi(t,y,y^\prime,y^{\prime\prime})=\alpha(t)\beta(y)\gamma(y^\prime)\delta(y^{\prime\prime})$, where $\alpha(t),\;\beta(y)\;\gamma(y^\prime),\;$ and $\delta(y^{\prime\prime})$ are differentiable functions. Substitute $\mu(\xi)=\mu(\alpha(t)\beta(y)\gamma(y^\prime)\delta(y^{\prime\prime}))$ in  \eqref{sys1}, then we get
\begin{equation}\label{sys2}
\left\{
\begin{array}{ll}
\mu(\xi)F_{3t}(\mathbf y)+\mu^\prime(\xi)\xi_t F_{3}(\mathbf y )=\mu(\xi )F_{0y^{\prime\prime}}(\mathbf y)+\mu^\prime(\xi )\xi_{y^{\prime\prime}}F_{0}(\mathbf y),\\
\mu(\xi)F_{2t}(\mathbf y)+\mu^\prime(\xi)\xi_tF_{2}(\mathbf y )=\mu(\xi )F_{0y^\prime}(\mathbf y)+\mu^\prime(\xi)\xi_{y^\prime}F_{0}(\mathbf y),\\
\mu(\xi)F_{1t}(\mathbf y)+\mu^\prime(\xi)\xi_tF_{1}(\mathbf y )=\mu(\xi )F_{0y}(\mathbf y)+\mu^\prime(\xi )\xi_yF_{0}(\mathbf y),\\
\mu(\xi)F_{1y^\prime}(\mathbf y)+\mu^\prime(\xi)\xi_{y^\prime} F_{1}(\mathbf y )=\mu(\xi)F_{2y}(\mathbf y)+\mu^\prime(\xi )\xi_yF_{2}(\mathbf y),\\
\mu(\xi)F_{1y^{\prime\prime}}(\mathbf y)+\mu^\prime(\xi) \xi_{y^{\prime\prime}}F_{1}(\mathbf y )=\mu(\xi )F_{2y}(\mathbf y)+\mu^\prime(\xi)\xi_yF_{2}(\mathbf y),\\
\mu(\xi)F_{2y^{\prime\prime}}(\mathbf y)+\mu^\prime(\xi)\xi_{y^{\prime\prime}}F_{2}(\mathbf y )=\mu(\xi)F_{3y^\prime}(\mathbf y)+\mu^\prime(\xi )\xi_{y^\prime}F_{3}(\mathbf y),
\end{array}
\right.
\end{equation}
where $\mu^\prime(\xi)=\frac{d\mu}{d\xi}.$
Equivalently, we have 
\begin{equation}\label{sys22}
\left\{
\begin{array}{ll}
\displaystyle \frac{\mu^\prime(\xi)}{\mu(\xi)}=\frac{F_{3t}(\mathbf y)-F_{0y^{\prime\prime}}(\mathbf y)}{\xi_{y^{\prime\prime}}F_{0}(\mathbf y)-\xi_t F_{3}(\mathbf y )},&
\displaystyle \frac{\mu^\prime(\xi)}{\mu(\xi)}=\frac{F_{2t}(\mathbf y)-F_{0y^\prime}(\mathbf y)}{\xi_{y^\prime}F_{0}(\mathbf y)-\xi_tF_{2}(\mathbf y )},\\\\
\displaystyle \frac{\mu^\prime(\xi)}{\mu(\xi)}=\frac{F_{1t}(\mathbf y)-F_{0y}(\mathbf y)}{\xi_yF_{0}(\mathbf y)-\xi_tF_{1}(\mathbf y )},&
\displaystyle \frac{\mu^\prime(\xi)}{\mu(\xi)}=\frac{F_{1y^\prime}(\mathbf y)-F_{2y}(\mathbf y)}{\xi_yF_{2}(\mathbf y)-\xi_{y^\prime} F_{1}(\mathbf y )},\\\\
\displaystyle \frac{\mu^\prime(\xi)}{\mu(\xi)}=\frac{F_{1y^{\prime\prime}}(\mathbf y)-F_{2y}(\mathbf y)}{\xi_yF_{2}(\mathbf y)- \xi_{y^{\prime\prime}}F_{1}(\mathbf y)} ,&
\displaystyle \frac{\mu^\prime(\xi)}{\mu(\xi)} =\frac{F_{2y^{\prime\prime}}(\mathbf y)-F_{3y^\prime}(\mathbf y)}{\xi_{y^{\prime\prime}}F_{2}(\mathbf y )-\xi_{y^\prime}F_{3}(\mathbf y)}.
\end{array}
\right.
\end{equation}
Hence, an integrating factor $\mu(\xi)$ of  \eqref{thirdorder} exists if 
\begin{itemize}
\item [a)] $\frac{F_{3t}(\mathbf y)-F_{0y^{\prime\prime}}(\mathbf y)}{\xi_{y^{\prime\prime}}F_{0}(\mathbf y)-\xi_t F_{3}(\mathbf y )}$, $\frac{F_{2t}(\mathbf y)-F_{0y^\prime}(\mathbf y)}{\xi_{y^\prime}F_{0}(\mathbf y)-\xi_tF_{2}(\mathbf y )}$, $\frac{F_{1t}(\mathbf y)-F_{0y}(\mathbf y)}{\xi_yF_{0}(\mathbf y)-\xi_tF_{1}(\mathbf y )}$, $\frac{F_{1y^\prime}(\mathbf y)-F_{2y}(\mathbf y)}{\xi_yF_{2}(\mathbf y)-\xi_{y^\prime} F_{1}(\mathbf y )}$, $\frac{F_{1y^{\prime\prime}}(\mathbf y)-F_{2y}(\mathbf y)}{\xi_yF_{2}(\mathbf y)- \xi_{y^{\prime\prime}}F_{1}(\mathbf y)}$, 

and $\frac{F_{2y^{\prime\prime}}(\mathbf y)-F_{3y^\prime}(\mathbf y)}{\xi_{y^{\prime\prime}}F_{2}(\mathbf y )-\xi_{y^\prime}F_{3}(\mathbf y)}$ are functions of $\xi,$

and
\item[b)] \begin{eqnarray*}
\frac{F_{3t}(\mathbf y)-F_{0y^{\prime\prime}}(\mathbf y)}{\xi_{y^{\prime\prime}}F_{0}(\mathbf y)-\xi_t F_{3}(\mathbf y )}&=& \frac{F_{2t}(\mathbf y)-F_{0y^\prime}(\mathbf y)}{\xi_{y^\prime}F_{0}(\mathbf y)-\xi_tF_{2}(\mathbf y )}= \frac{F_{1t}(\mathbf y)-F_{0y}(\mathbf y)}{\xi_yF_{0}(\mathbf y)-\xi_tF_{1}(\mathbf y )}= \frac{F_{1y^\prime}(\mathbf y)-F_{2y}(\mathbf y)}{\xi_yF_{2}(\mathbf y)-\xi_{y^\prime} F_{1}(\mathbf y )}=\\ \frac{F_{1y^{\prime\prime}}(\mathbf y)-F_{2y}(\mathbf y)}{\xi_yF_{2}(\mathbf y)- \xi_{y^{\prime\prime}}F_{1}(\mathbf y)}&=&\frac{F_{2y^{\prime\prime}}(\mathbf y)-F_{3y^\prime}(\mathbf y)}{\xi_{y^{\prime\prime}}F_{2}(\mathbf y )-\xi_{y^\prime}F_{3}(\mathbf y)}
\end{eqnarray*}
\end{itemize}
hold. Therefore, we have the following theorem:
\begin{theorem}\label{maintheorem} Let $\mathbf{y}=(t,y,y^\prime,y^{\prime\prime})$, and 
assume that Equation \eqref{thirdorder} is none exact differential equation. Then it  admits an integrating factor $\mu(\xi)=\mu(\xi(t,y,y^\prime,y^{\prime\prime}))=\mu\left(\alpha(t)\beta(y)\gamma(y^\prime)\delta(y^{\prime\prime})\right)$, where $\alpha(t), \beta(y), \gamma(y^\prime),$ and  $\delta(y^{\prime\prime})$ are differentiable functions; if the  conditions
\begin{itemize}
\item [a)] $\frac{F_{3t}(\mathbf y)-F_{0y^{\prime\prime}}(\mathbf y)}{\xi_{y^{\prime\prime}}F_{0}(\mathbf y)-\xi_t F_{3}(\mathbf y )}$, $\frac{F_{2t}(\mathbf y)-F_{0y^\prime}(\mathbf y)}{\xi_{y^\prime}F_{0}(\mathbf y)-\xi_tF_{2}(\mathbf y )}$, $\frac{F_{1t}(\mathbf y)-F_{0y}(\mathbf y)}{\xi_yF_{0}(\mathbf y)-\xi_tF_{1}(\mathbf y )}$, $\frac{F_{1y^\prime}(\mathbf y)-F_{2y}(\mathbf y)}{\xi_yF_{2}(\mathbf y)-\xi_{y^\prime} F_{1}(\mathbf y )}$, $\frac{F_{1y^{\prime\prime}}(\mathbf y)-F_{2y}(\mathbf y)}{\xi_yF_{2}(\mathbf y)- \xi_{y^{\prime\prime}}F_{1}(\mathbf y)}$, 

and $\frac{F_{2y^{\prime\prime}}(\mathbf y)-F_{3y^\prime}(\mathbf y)}{\xi_{y^{\prime\prime}}F_{2}(\mathbf y )-\xi_{y^\prime}F_{3}(\mathbf y)}$ are functions in $\xi:=\alpha(t)\beta(y)\gamma(y^\prime)\delta(y^{\prime\prime}),$

and 
\item[b)] \begin{eqnarray*}
\frac{F_{3t}(\mathbf y)-F_{0y^{\prime\prime}}(\mathbf y)}{\xi_{y^{\prime\prime}}F_{0}(\mathbf y)-\xi_t F_{3}(\mathbf y )}&=& \frac{F_{2t}(\mathbf y)-F_{0y^\prime}(\mathbf y)}{\xi_{y^\prime}F_{0}(\mathbf y)-\xi_tF_{2}(\mathbf y )}= \frac{F_{1t}(\mathbf y)-F_{0y}(\mathbf y)}{\xi_yF_{0}(\mathbf y)-\xi_tF_{1}(\mathbf y )}= \frac{F_{1y^\prime}(\mathbf y)-F_{2y}(\mathbf y)}{\xi_yF_{2}(\mathbf y)-\xi_{y^\prime} F_{1}(\mathbf y )}=\\ \frac{F_{1y^{\prime\prime}}(\mathbf y)-F_{2y}(\mathbf y)}{\xi_yF_{2}(\mathbf y)- \xi_{y^{\prime\prime}}F_{1}(\mathbf y)}&=&\frac{F_{2y^{\prime\prime}}(\mathbf y)-F_{3y^\prime}(\mathbf y)}{\xi_{y^{\prime\prime}}F_{2}(\mathbf y )-\xi_{y^\prime}F_{3}(\mathbf y)},
\end{eqnarray*}
\end{itemize}
hold.
Moreover, the integrating factor is given by the formula
\[
\mu(\xi)=\exp\left\{\displaystyle\int\frac{F_{3t}(\mathbf y)-F_{0y^{\prime\prime}}(\mathbf y)}{\xi_{y^{\prime\prime}}F_{0}(\mathbf y)-\xi_t F_{3}(\mathbf y )}d\xi\right\}.
\]
\end{theorem}
In the following sections, we presents some special cases of the above theorem. Moreover, we present some examples to illustrate the methodology of finding such  integrating factor.
\subsection{Integrating Factors of the Forms $\mu(\alpha(t)), \mu(\beta(y)), \mu(\gamma(y^\prime))$ and $\mu(\delta(y^{\prime\prime}))$}
In this section, we give conditions so that an integrating factor of one of the forms  $\mu(\alpha(t))$, $\mu(\beta(y))$, $\mu(\gamma(y^\prime))$ and $\mu(\delta(y^{\prime\prime}))$ for equation \eqref{thirdorder} exists. As a result of Theorem \ref{maintheorem}, we have  the following corollaries: 

\begin{corollary}\label{cort}
 Let $\mathbf{y}=(t,y,y^\prime,y^{\prime\prime})$, and 
assume that Equation \eqref{thirdorder} is none exact differential equation. Then it  admits an integrating factor $\mu(\xi )=\mu(\alpha(t))$, where $\alpha(t)$ is  differentiable function; if the following two conditions hold:
\begin{itemize}
\item [a)] 
$
F_{1y^\prime}(\mathbf y)=F_{2y}(\mathbf y),$ $F_{1y^{\prime\prime}}(\mathbf y)=F_{2y}(\mathbf y),$ and 
$F_{2y^{\prime\prime}}(\mathbf y)=F_{3y^\prime}(\mathbf y),$

and 
\item[b)] 
$\frac{F_{0y^{\prime\prime}}(\mathbf y)-F_{3t}(\mathbf y)}{\xi_t F_{3}(\mathbf y )}= \frac{F_{0y^\prime}(\mathbf y)-F_{2t}(\mathbf y)}{\xi_tF_{2}(\mathbf y )}=  
\frac{F_{0y}(\mathbf y)-F_{1t}(\mathbf y)}{\xi_tF_{1}(\mathbf y )}$, and they are functions in $\xi:=\xi(\alpha(t)).$ 
\end{itemize}
Moreover, the integrating factor is given by 
\[
\mu(\xi)=\exp\left\{\displaystyle\int\frac{F_{0y^{\prime\prime}}(\mathbf y)-F_{3t}(\mathbf y)}{\xi_t F_{3}(\mathbf y )}d\xi\right\}.
\]
\end{corollary}

\begin{corollary}
 Let $\mathbf{y}=(t,y,y^\prime,y^{\prime\prime})$, and 
assume that Equation \eqref{thirdorder} is none exact differential equation. Then it  admits an integrating factor $\mu(\xi )=\mu(\beta(y))$, where $\beta(y)$ is  differentiable function; if the following two conditions hold:
\begin{itemize}
\item [a)] $F_{3t}(\mathbf y)=F_{0y^{\prime\prime}}(\mathbf y)$, $F_{2t}(\mathbf y)=F_{0y^\prime}(\mathbf y)$, and  $F_{2y^{\prime\prime}}(\mathbf y)=F_{3y^\prime}(\mathbf y),$ 

and 
\item[b)] 
$\frac{F_{1t}(\mathbf y)-F_{0y}(\mathbf y)}{\xi_yF_{0}(\mathbf y)}= \frac{F_{1y^\prime}(\mathbf y)-F_{2y}(\mathbf y)}{\xi_yF_{2}(\mathbf y)}= \frac{F_{1y^{\prime\prime}}(\mathbf y)-F_{2y}(\mathbf y)}{\xi_yF_{2}(\mathbf y)},$ and they are functions in $\xi:=\xi(\beta(y)).$
\end{itemize}
Moreover, the integrating factor is given by 
\[
\mu(\xi)=\exp\left\{\displaystyle\int\frac{F_{1t}(\mathbf y)-F_{0y}(\mathbf y)}{\xi_yF_{0}(\mathbf y)}d\xi\right\}.
\]
\end{corollary}

\begin{corollary}\label{coryprime}
 Let $\mathbf{y}=(t,y,y^\prime,y^{\prime\prime})$, and 
assume that Equation \eqref{thirdorder} is none exact differential equation. Then it  admits an integrating factor $\mu(\xi)=\mu(\gamma(y^\prime))$ for some differentiable functions $\gamma(y^\prime)$; if the two conditions hold:
\begin{itemize}
\item [a)] $F_{3t}(\mathbf y)=F_{0y^{\prime\prime}}(\mathbf y)$,  $F_{1t}(\mathbf y)=F_{0y}(\mathbf y)$,  and $F_{1y^{\prime\prime}}(\mathbf y)=F_{2y}(\mathbf y)$, 

and 
\item[b)] $
\frac{F_{2t}(\mathbf y)-F_{0y^\prime}(\mathbf y)}{\xi_{y^\prime}F_{0}(\mathbf y)}= \frac{F_{2y}(\mathbf y)-F_{1y^\prime}(\mathbf y)}{\xi_{y^\prime} F_{1}(\mathbf y )}=\frac{F_{3y^\prime}(\mathbf y)-F_{2y^{\prime\prime}}(\mathbf y)}{\xi_{y^\prime}F_{3}(\mathbf y)},$ and they are functions in $\xi:=\gamma(y^\prime)$.
\end{itemize}
Moreover, the integrating factor is given by 
\[
\mu(\xi)=\exp\left\{\displaystyle\int\frac{F_{2t}(\mathbf y)-F_{0y^\prime}(\mathbf y)}{\xi_{y^\prime}F_{0}(\mathbf y)}d\xi\right\}.
\]
\end{corollary}

\begin{corollary}
 Let $\mathbf{y}=(t,y,y^\prime,y^{\prime\prime})$, and 
assume that Equation \eqref{thirdorder} is none exact differential equation. Then it  admits an integrating factor $\mu(\xi )=\mu(\delta(y^{\prime\prime}))$ for some differentiable function  $\delta(y^{\prime\prime})$; if the following two conditions hold:
\begin{itemize}
\item [a)]  $F_{2t}(\mathbf y)=F_{0y^\prime}(\mathbf y)$, $F_{1t}(\mathbf y)=F_{0y}(\mathbf y)$, and $F_{1y^\prime}(\mathbf y)=F_{2y}(\mathbf y)$, 

and 
\item[b)] $
\frac{F_{3t}(\mathbf y)-F_{0y^{\prime\prime}}(\mathbf y)}{\xi_{y^{\prime\prime}}F_{0}(\mathbf y)}= \frac{F_{2y}(\mathbf y)-F_{1y^{\prime\prime}}(\mathbf y)}{ \xi_{y^{\prime\prime}}F_{1}(\mathbf y)}=\frac{F_{2y^{\prime\prime}}(\mathbf y)-F_{3y^\prime}(\mathbf y)}{\xi_{y^{\prime\prime}}F_{2}(\mathbf y )},
$ and they are functions in $\xi:=\delta(y^{\prime\prime}).$
\end{itemize}
Moreover, the integrating factor is given by 
\[
\mu(\xi)=\exp\left\{\displaystyle\int\frac{F_{3t}(\mathbf y)-F_{0y^{\prime\prime}}(\mathbf y)}{\xi_{y^{\prime\prime}}F_{0}(\mathbf y)-\xi_t F_{3}(\mathbf y )}d\xi\right\}.
\]
\end{corollary}

\begin{example}
Consider the differential equation 
\begin{equation}\label{exampleeq}
(y^\prime)^3y^{\prime\prime\prime}+2yy^{\prime\prime}-(y^\prime)^2+(y^\prime)^3=0.
\end{equation}
Clearly, 
$F_3(t,y,y^\prime,y^{\prime\prime})=(y^\prime)^3$, $F_2(t,y,y^\prime,y^{\prime\prime})=2y$, $F_1(t,y,y^\prime,y^{\prime\prime})=-y^\prime$, and $F_0(t,y,y^\prime,y^{\prime\prime})=(y^\prime)^3$. Moreover,
$F_{3t}(t,y,y^\prime,y^{\prime\prime})=F_{0y^{\prime\prime}}(t,y,y^\prime,y^{\prime\prime})=0$,  $F_{1t}(t,y,y^\prime,y^{\prime\prime})=F_{0y}(t,y,y^\prime,y^{\prime\prime})=0$, and 
$
F_{1y^{\prime\prime}}(t,y,y^\prime,y^{\prime\prime})=F_{2y}(t,y,y^\prime,y^{\prime\prime})=0$, and 
$
\frac{F_{2t}(\mathbf y)-F_{0y^\prime}(\mathbf y)}{\xi_{y^\prime}F_{0}(\mathbf y)}= \frac{F_{2y}(\mathbf y)-F_{1y^\prime}(\mathbf y)}{\xi_{y^\prime} F_{1}(\mathbf y )}=\frac{F_{3y^\prime}(\mathbf y)-F_{2y^{\prime\prime}}(\mathbf y)}{\xi_{y^\prime}F_{3}(\mathbf y)}=\frac{-3}{y^\prime}.$ Therefor, the condition in Corollary \ref{coryprime} hold, and so an integrating factor of this third order differential equation exists and is given by the formula $\mu(y^\prime)=\left(y^\prime\right)^{-3}$. By multiplying \eqref{exampleeq} by
$\mu(y^\prime)=\left(y^\prime\right)^{-3}$, we get 
\begin{equation}
y^{\prime\prime\prime}+2y\left(y^\prime\right)^{-3}y^{\prime\prime}-(y^\prime)^{-2}y^\prime+1=0.
\end{equation}
For this equation $F_3(t,y,y^\prime,y^{\prime\prime})=1$, $F_2(t,y,y^\prime,y^{\prime\prime})=2y\left(y^\prime\right)^{-3}$, $F_1(t,y,y^\prime,y^{\prime\prime})=-(y^\prime)^{-2}$, and $F_0(t,y,y^\prime,y^{\prime\prime})=1.$ Clearly, $\partial_{y^{\prime\prime}}F_0=\partial_tF_3=0,\;\partial_{y^{\prime\prime}}
F_1=\partial_yF_3=0,\; \partial_{y^{\prime\prime}}F_2=\partial_{y^\prime}F_3=0,
\partial_{y^\prime}F_0=\partial_tF_2=0,\;\partial_{y^\prime}F_1=\partial_yF_2=2(y^\prime)^{-3},\; 
\text{and}\;\partial_yF_0=\partial_tF_1=0.$ Hence, it is exact differential equation, and its first integral is given by
\begin{equation}
 \Psi\left(t,y,y^\prime,y^{\prime\prime}\right)= \int_{t_0}^{t}d\xi-(y^\prime)^{-2}\int_{y_0}^{y} d\xi
+2y_0\int_{y^\prime_0}^{y^\prime}\xi^{-3}d\xi +\int_{y_0^{\prime\prime}}^{y^{\prime\prime}} d\xi 
=c,
\end{equation}
More precisely, 
 \begin{equation}
 \Psi\left(t,y,y^\prime,y^{\prime\prime}\right)=y^{\prime\prime}-(y^\prime)^{-2}y+t=c
%
\end{equation}
\end{example}

\begin{example}
consider the third order linear differential equation
\begin{equation}\label{remlineareq}
p_2(t)y^{\prime\prime\prime}+\alpha p_2(t)y^{\prime\prime}+p_1(t)y^\prime+p_0(t)y=h(t),\;p_2(t)\neq 0.
\end{equation}
Assume that $p_2(t)p_1^\prime(t)-p_2^\prime(t) p_1(t)=p_0(t)p_2(t)$. Then, by applying the conditions in Corollary \ref{cort}, the above third order linear differential equation admits an integrating factor $\mu(t)=\frac{1}{p_2(t)}$. In fact, if we multiply Eq. \eqref{remlineareq} by $\mu(t)=\frac{1}{p_2(t)}$, then we get
\[
y^{\prime\prime\prime}+\alpha y^{\prime\prime}+\frac{p_1(t)}{p_2(t)}y^\prime+\frac{p_0(t)}{p_2(t)}y=\frac{h(t)}{p_2(t)}.
\]
From the condition $p_2(t)p_1^\prime(t)-p_2^\prime(t) p_1(t)=p_0(t)p_2(t)$, we have $\left(\frac{p_1(t)}{p_2(t)}\right)^\prime=\frac{p_0(t)}{p_2(t)}$. Hence, the above equation becomes
\[
y^{\prime\prime\prime}+\alpha y^{\prime\prime}+\left(\frac{p_1(t)}{p_2(t)}\right)y^\prime+\left(\frac{p_1(t)}{p_2(t)}\right)^\prime y=\frac{h(t)}{p_2(t)}.
\] 
This equation can be written as
\[
\frac{d}{dt}\left[y^{\prime\prime}+\alpha y^{\prime}+\left(\frac{p_1(t)}{p_2(t)}\right)y\right]=\frac{h(t)}{p_2(t)}.
\]
Hence, the first integral of Eq. \eqref{remlineareq} is given by
\[
y^{\prime\prime}+\alpha y^{\prime}+\left(\frac{p_1(t)}{p_2(t)}\right)y=\int^t\frac{h(s)}{p_2(s)ds}+c.
\]
\end{example}

\section{Integrating Factors of the Forms $\mu(\alpha(t)\beta(y)),$ $\mu(\alpha(t)\gamma(y^\prime)),$ $\mu(\alpha(t)\delta(y^{\prime\prime})),$ $\mu(\beta(y)\gamma(y^\prime)),$ $\mu(\beta(y)\delta(y^{\prime\prime}))$ and $\mu(\gamma(y^\prime)\delta(y^{\prime\prime}))$}

In this section, we give conditions so that an integrating factor of one of the forms  $\mu(\alpha(t)\beta(y)),$ $\mu(\alpha(t)\gamma(y^\prime)),$ $\mu(\alpha(t)\delta(y^{\prime\prime})),$ $\mu(\beta(y)\gamma(y^\prime)),$ $\mu(\beta(y)\delta(y^{\prime\prime}))$ and $\mu(\gamma(y^\prime)\delta(y^{\prime\prime}))$ for equation \eqref{thirdorder} exists. As a result of Theorem \ref{maintheorem}, we have  the following corollaries: 

\begin{corollary}
 Let $\mathbf{y}=(t,y,y^\prime,y^{\prime\prime})$, and 
assume that Equation \eqref{thirdorder} is none exact differential equation. Then it  admits an integrating factor $\mu(\xi )=\mu(\alpha(t)\beta(y))$, where $\alpha(t)$ and $\beta(y))$ are  differentiable functions; if the following two conditions hold:
\begin{itemize}
\item [a)] 
$F_{2y^{\prime\prime}}(\mathbf y)=F_{3y^\prime}(\mathbf y),$
and 
\item[b)] 
$
\frac{F_{0y^{\prime\prime}}(\mathbf y)-F_{3t}(\mathbf y)}{\xi_t F_{3}(\mathbf y )}= \frac{F_{0y^\prime}(\mathbf y)-F_{2t}(\mathbf y)}{\xi_tF_{2}(\mathbf y )}= \frac{F_{1y^\prime}(\mathbf y)-F_{2y}(\mathbf y)}{\xi_yF_{2}(\mathbf y)}= \frac{F_{1y^{\prime\prime}}(\mathbf y)-F_{2y}(\mathbf y)}{\xi_yF_{2}(\mathbf y)}= \frac{F_{1t}(\mathbf y)-F_{0y}(\mathbf y)}{\xi_yF_{0}(\mathbf y)-\xi_tF_{1}(\mathbf y )}$, and they are functions in $\xi:=\xi(\alpha(t)\beta(y)).$ 
\end{itemize}
Moreover, the integrating factor is given by 
\[
\mu(\xi)=\exp\left\{\displaystyle\int\frac{F_{0y^{\prime\prime}}(\mathbf y)-F_{3t}(\mathbf y)}{\xi_t F_{3}(\mathbf y )}d\xi\right\}.
\]
\end{corollary}

\begin{corollary}
 Let $\mathbf{y}=(t,y,y^\prime,y^{\prime\prime})$, and 
assume that Equation \eqref{thirdorder} is none exact differential equation. Then it  admits an integrating factor $\mu(\xi)=\mu(\alpha(t)\gamma(y^\prime))$ for some differentiable functions $\alpha(t)$ and $\gamma(y^\prime)$; if the two conditions hold:
\begin{itemize}
\item [a)] $F_{1y^{\prime\prime}}(\mathbf y)=F_{2y}(\mathbf y)$, 

and 
\item[b)] $
\frac{F_{0y^{\prime\prime}}(\mathbf y)-F_{3t}(\mathbf y)}{\xi_t F_{3}(\mathbf y )}=
\frac{F_{0y}(\mathbf y)-F_{1t}(\mathbf y)}{\xi_tF_{1}(\mathbf y )}= \frac{F_{2y}(\mathbf y)-F_{1y^\prime}(\mathbf y)}{\xi_{y^\prime} F_{1}(\mathbf y )}=
\frac{F_{3y^\prime}(\mathbf y)-F_{2y^{\prime\prime}}(\mathbf y)}{\xi_{y^\prime}F_{3}(\mathbf y)}=
\frac{F_{2t}(\mathbf y)-F_{0y^\prime}(\mathbf y)}{\xi_{y^\prime}F_{0}(\mathbf y)-\xi_tF_{2}(\mathbf y )}.
$ and they are functions in $\xi:=\alpha(t)\gamma(y^\prime)$.
\end{itemize}
Moreover, the integrating factor is given by 
\[
\mu(\xi)=\exp\left\{\displaystyle\int\frac{F_{0y^{\prime\prime}}(\mathbf y)-F_{3t}(\mathbf y)}{\xi_t F_{3}(\mathbf y )}d\xi\right\}.
\]
\end{corollary}

\begin{corollary}
 Let $\mathbf{y}=(t,y,y^\prime,y^{\prime\prime})$, and 
assume that Equation \eqref{thirdorder} is none exact differential equation. Then it  admits an integrating factor $\mu(\xi )=\mu(\alpha(t)\delta(y^{\prime\prime}))$ for some differentiable function  $\delta(y^{\prime\prime})$; if the following two conditions hold:
\begin{itemize}
\item [a)]   $F_{1y^\prime}(\mathbf y)=F_{2y}(\mathbf y)$, 

and 
\item[b)] $
\frac{F_{0y^\prime}(\mathbf y)-F_{2t}(\mathbf y)}{\xi_tF_{2}(\mathbf y )}=
\frac{F_{0y}(\mathbf y)-F_{1t}(\mathbf y)}{\xi_tF_{1}(\mathbf y )}= 
\frac{F_{2y}(\mathbf y)-F_{1y^{\prime\prime}}(\mathbf y)}{ \xi_{y^{\prime\prime}}F_{1}(\mathbf y)}=
\frac{F_{2y^{\prime\prime}}(\mathbf y)-F_{3y^\prime}(\mathbf y)}{\xi_{y^{\prime\prime}}F_{2}(\mathbf y )}=
\frac{F_{3t}(\mathbf y)-F_{0y^{\prime\prime}}(\mathbf y)}{\xi_{y^{\prime\prime}}F_{0}(\mathbf y)}.
$ and they are functions in $\xi:=\alpha(t)\delta(y^{\prime\prime}).$
\end{itemize}
Moreover, the integrating factor is given by 
\[
\mu(\xi)=\exp\left\{\displaystyle\int\frac{F_{0y^\prime}(\mathbf y)-F_{2t}(\mathbf y)}{\xi_tF_{2}(\mathbf y )}d\xi\right\}.
\]
\end{corollary}

\begin{corollary}
 Let $\mathbf{y}=(t,y,y^\prime,y^{\prime\prime})$, and 
assume that Equation \eqref{thirdorder} is none exact differential equation. Then it  admits an integrating factor $\mu(\xi )=\mu(\beta(y)\gamma(y^\prime))$, where $\beta(y)$ and $\gamma(y^\prime)$ are  differentiable functions; if the following two conditions hold:
\begin{itemize}
\item [a)] $F_{3t}(\mathbf y)=F_{0y^{\prime\prime}}$
and 
\item[b)] 
$
\frac{F_{2t}(\mathbf y)-F_{0y^\prime}(\mathbf y)}{\xi_{y^\prime}F_{0}(\mathbf y)}= 
\frac{F_{1t}(\mathbf y)-F_{0y}(\mathbf y)}{\xi_yF_{0}(\mathbf y)}=
\frac{F_{1y^{\prime\prime}}(\mathbf y)-F_{2y}(\mathbf y)}{\xi_yF_{2}(\mathbf y)}=
\frac{F_{3y^\prime}(\mathbf y)-F_{2y^{\prime\prime}}(\mathbf y)}{\xi_{y^\prime}F_{3}(\mathbf y)}=
\frac{F_{1y^\prime}(\mathbf y)-F_{2y}(\mathbf y)}{\xi_yF_{2}(\mathbf y)-\xi_{y^\prime} F_{1}(\mathbf y )}
$ and they are functions in $\xi:=\beta(y)\gamma(y^\prime).$
\end{itemize}
Moreover, the integrating factor is given by 
\[
\mu(\xi)=\exp\left\{\displaystyle\int\frac{F_{2t}(\mathbf y)-F_{0y^\prime}(\mathbf y)}{\xi_{y^\prime}F_{0}(\mathbf y)}d\xi\right\}.
\]
\end{corollary}

\begin{corollary}
 Let $\mathbf{y}=(t,y,y^\prime,y^{\prime\prime})$, and 
assume that Equation \eqref{thirdorder} is none exact differential equation. Then it  admits an integrating factor $\mu(\xi )=\mu(\beta(y)\delta(y^{\prime\prime}))$, where $\beta(y)$ and $\delta(y^{\prime\prime})$ are  differentiable functions; if the following two conditions hold:
\begin{itemize}
\item [a)]  $F_{2t}(\mathbf y)=F_{0y^\prime}(\mathbf y)$,
and 
\item[b)] 
$
\frac{F_{3t}(\mathbf y)-F_{0y^{\prime\prime}}(\mathbf y)}{\xi_{y^{\prime\prime}}F_{0}(\mathbf y)}=
\frac{F_{1t}(\mathbf y)-F_{0y}(\mathbf y)}{\xi_yF_{0}(\mathbf y)}= 
\frac{F_{1y^\prime}(\mathbf y)-F_{2y}(\mathbf y)}{\xi_yF_{2}(\mathbf y)}=
\frac{F_{2y^{\prime\prime}}(\mathbf y)-F_{3y^\prime}(\mathbf y)}{\xi_{y^{\prime\prime}}F_{2}(\mathbf y )}= 
\frac{F_{1y^{\prime\prime}}(\mathbf y)-F_{2y}(\mathbf y)}{\xi_yF_{2}(\mathbf y)- \xi_{y^{\prime\prime}}F_{1}(\mathbf y)}
$ and they are functions in $\xi:=\beta(y)\delta(y^{\prime\prime}).$
\end{itemize}
Moreover, the integrating factor is given by 
\[
\mu(\xi)=\exp\left\{\displaystyle\int\frac{F_{3t}(\mathbf y)-F_{0y^{\prime\prime}}(\mathbf y)}{\xi_{y^{\prime\prime}}F_{0}(\mathbf y)}d\xi\right\}.
\]
\end{corollary}

\begin{corollary}
Let $\mathbf{y}=(t,y,y^\prime,y^{\prime\prime})$, and 
assume that Equation \eqref{thirdorder} is none exact differential equation. Then it  admits an integrating factor $\mu(\xi )=\mu(\gamma(y^\prime)\delta(y^{\prime\prime}))$, where $\gamma(y^\prime)$ and $\delta(y^{\prime\prime})$ are  differentiable functions; if the following two conditions hold:
\begin{itemize}
\item [a)] $F_{1t}(\mathbf y)=F_{0y}(\mathbf y)$, 
and 
\item[b)] 
$
\frac{F_{3t}(\mathbf y)-F_{0y^{\prime\prime}}(\mathbf y)}{\xi_{y^{\prime\prime}}F_{0}(\mathbf y)}= 
\frac{F_{2t}(\mathbf y)-F_{0y^\prime}(\mathbf y)}{\xi_{y^\prime}F_{0}(\mathbf y)}= 
\frac{F_{2y}(\mathbf y)-F_{1y^\prime}(\mathbf y)}{\xi_{y^\prime} F_{1}(\mathbf y )}=
\frac{F_{2y}(\mathbf y)-F_{1y^{\prime\prime}}(\mathbf y)}{ \xi_{y^{\prime\prime}}F_{1}(\mathbf y)}=
\frac{F_{2y^{\prime\prime}}(\mathbf y)-F_{3y^\prime}(\mathbf y)}{\xi_{y^{\prime\prime}}F_{2}(\mathbf y )-\xi_{y^\prime}F_{3}(\mathbf y)}
$ and they are functions in $\xi=\gamma(y^\prime)\delta(y^{\prime\prime}).$
\end{itemize}
Moreover, the integrating factor is given by 
\[
\mu(\xi)=\exp\left\{\displaystyle\int\frac{F_{3t}(\mathbf y)-F_{0y^{\prime\prime}}(\mathbf y)}{\xi_{y^{\prime\prime}}F_{0}(\mathbf y)}d\xi\right\}.
\]
\end{corollary}
\section{Integrating Factors of the Forms $\mu(\alpha(t)\beta(y)\gamma(y^\prime)),$ $\mu(\alpha(t)\beta(y)\delta(y^{\prime\prime})),$  $\mu(\alpha(t)\gamma(y^\prime)\delta(y^{\prime\prime}))$ and $\mu(\beta(y)\gamma(y^\prime)\delta(y^{\prime\prime}))$}

In this section, we give conditions so that an integrating factor of one of the forms  $\mu(\alpha(t)\beta(y)\gamma(y^\prime)),$ $\mu(\alpha(t)\beta(y)\delta(y^{\prime\prime})),$  $\mu(\alpha(t)\gamma(y^\prime)\delta(y^{\prime\prime}))$ and $\mu(\beta(y)\gamma(y^\prime)\delta(y^{\prime\prime}))$
 for equation \eqref{thirdorder} exists. As a result of Theorem \ref{maintheorem}, we have  the following corollaries: 

\begin{corollary}
 Let $\mathbf{y}=(t,y,y^\prime,y^{\prime\prime})$, and 
assume that Equation \eqref{thirdorder} is none exact differential equation. Then it  admits an integrating factor $\mu(\xi)=\mu(\xi(t,y,y^\prime))=\mu\left(\alpha(t)\beta(y)\gamma(y^\prime)\right)$ where $\alpha(t), \beta(y),$ and $ \gamma(y^\prime)$ are differentiable functions; if
\begin{eqnarray*} 
\frac{F_{0y^{\prime\prime}}(\mathbf y)-F_{3t}(\mathbf y)}{\xi_t F_{3}(\mathbf y )}&=&\frac{F_{1y^{\prime\prime}}(\mathbf y)-F_{2y}(\mathbf y)}{\xi_yF_{2}(\mathbf y)}=\frac{F_{3y^\prime}(\mathbf y)-F_{2y^{\prime\prime}}(\mathbf y)}{\xi_{y^\prime}F_{3}(\mathbf y)}=
\frac{F_{2t}(\mathbf y)-F_{0y^\prime}(\mathbf y)}{\xi_{y^\prime}F_{0}(\mathbf y)-\xi_tF_{2}(\mathbf y )}=\\\frac{F_{1t}(\mathbf y)-F_{0y}(\mathbf y)}{\xi_yF_{0}(\mathbf y)-\xi_tF_{1}(\mathbf y )}&=&\frac{F_{1y^\prime}(\mathbf y)-F_{2y}(\mathbf y)}{\xi_yF_{2}(\mathbf y)-\xi_{y^\prime} F_{1}(\mathbf y )}
\end{eqnarray*}
and they are functions in $\xi:=\alpha(t)\beta(y)\gamma(y^\prime).$ Moreover, the integrating factor is given by the formula
\[
\mu(\xi)=\exp\left\{\displaystyle\int\frac{F_{0y^{\prime\prime}}(\mathbf y)-F_{3t}(\mathbf y)}{\xi_t F_{3}(\mathbf y )}d\xi\right\}.
\]
\end{corollary}

\begin{corollary}
 Let $\mathbf{y}=(t,y,y^\prime,y^{\prime\prime})$, and 
assume that Equation \eqref{thirdorder} is none exact differential equation. Then it  admits an integrating factor $\mu(\xi)=\mu(\alpha(t)\beta(y)\delta(y^{\prime\prime}))$ where $\alpha(t)$, $\beta(y)$, and $\delta(y^{\prime\prime})$ are differentiable functions; if 
\begin{eqnarray*}
\frac{F_{0y^\prime}(\mathbf y)-F_{2t}(\mathbf y)}{\xi_tF_{2}(\mathbf y )}
&=& 
\frac{F_{2y^{\prime\prime}}(\mathbf y)-F_{3y^\prime}(\mathbf y)}{\xi_{y^{\prime\prime}}F_{2}(\mathbf y )} = 
\frac{F_{1y^\prime}(\mathbf y)-F_{2y}(\mathbf y)}{\xi_yF_{2}(\mathbf y)}=
\frac{F_{1y^{\prime\prime}}(\mathbf y)-F_{2y}(\mathbf y)}{\xi_yF_{2}(\mathbf y)- \xi_{y^{\prime\prime}}F_{1}(\mathbf y)}=\\\frac{F_{3t}(\mathbf y)-F_{0y^{\prime\prime}}(\mathbf y)}{\xi_{y^{\prime\prime}}F_{0}(\mathbf y)-\xi_t F_{3}(\mathbf y )}&=&\frac{F_{1t}(\mathbf y)-F_{0y}(\mathbf y)}{\xi_yF_{0}(\mathbf y)-\xi_tF_{1}(\mathbf y )}
\end{eqnarray*}
and they are functions in $\xi=\alpha(t)\beta(y)\delta(y^{\prime\prime}).$
Moreover, the integrating factor is given by 
\[
\mu(\xi)=\exp\left\{\displaystyle\int\frac{F_{0y^\prime}(\mathbf y)-F_{2t}(\mathbf y)}{\xi_tF_{2}(\mathbf y )}d\xi\right\}.
\]
\end{corollary}

\begin{corollary}
 Let $\mathbf{y}=(t,y,y^\prime,y^{\prime\prime})$, and 
assume that Equation \eqref{thirdorder} is none exact differential equation. Then it  admits an integrating factor $\mu(\xi )=\mu(\alpha(t)\gamma(y^\prime)\delta(y^{\prime\prime}))$ where $\alpha(t)$, $\gamma(y^\prime)$, and $\delta(y^{\prime\prime})$ are differentiable functions; if 
\begin{eqnarray*}
\frac{F_{0y}(\mathbf y)-F_{1t}(\mathbf y)}{\xi_tF_{1}(\mathbf y )}&=& 
\frac{F_{2y}(\mathbf y)-F_{1y^\prime}(\mathbf y)}{\xi_{y^\prime} F_{1}(\mathbf y )}= 
\frac{F_{2y}(\mathbf y)-F_{1y^{\prime\prime}}(\mathbf y)}{ \xi_{y^{\prime\prime}}F_{1}(\mathbf y)}
\frac{F_{3t}(\mathbf y)-F_{0y^{\prime\prime}}(\mathbf y)}{\xi_{y^{\prime\prime}}F_{0}(\mathbf y)-\xi_t F_{3}(\mathbf y )}=\\ \frac{F_{2t}(\mathbf y)-F_{0y^\prime}(\mathbf y)}{\xi_{y^\prime}F_{0}(\mathbf y)-\xi_tF_{2}(\mathbf y )} 
&=&
\frac{F_{2y^{\prime\prime}}(\mathbf y)-F_{3y^\prime}(\mathbf y)}{\xi_{y^{\prime\prime}}F_{2}(\mathbf y )-\xi_{y^\prime}F_{3}(\mathbf y)}
\end{eqnarray*} 
and they are functions in $\xi=\alpha(t)\gamma(y^\prime)\delta(y^{\prime\prime}).$ Moreover, the integrating factor is given by 
\[
\mu(\xi)=\exp\left\{\displaystyle\int
\frac{F_{0y}(\mathbf y)-F_{1t}(\mathbf y)}{\xi_tF_{1}(\mathbf y )}d\xi\right\}.
\]
\end{corollary}

\begin{corollary}
 Let $\mathbf{y}=(t,y,y^\prime,y^{\prime\prime})$, and 
assume that Equation \eqref{thirdorder} is none exact differential equation. Then it  admits an integrating factor $\mu(\xi )=\mu(\beta(y)\gamma(y^\prime)\delta(y^{\prime\prime}))$, where $\beta(y)$, $\gamma(y^\prime)$ and $\delta(y^{\prime\prime})$ are differentiable functions; if 
 
\begin{eqnarray*}
\frac{F_{3t}(\mathbf y)-F_{0y^{\prime\prime}}(\mathbf y)}{\xi_{y^{\prime\prime}}F_{0}(\mathbf y)}&=& \frac{F_{2t}(\mathbf y)-F_{0y^\prime}(\mathbf y)}{\xi_{y^\prime}F_{0}(\mathbf y)}= \frac{F_{1t}(\mathbf y)-F_{0y}(\mathbf y)}{\xi_yF_{0}(\mathbf y)}= \frac{F_{1y^\prime}(\mathbf y)-F_{2y}(\mathbf y)}{\xi_yF_{2}(\mathbf y)-\xi_{y^\prime} F_{1}(\mathbf y )}=\\ \frac{F_{1y^{\prime\prime}}(\mathbf y)-F_{2y}(\mathbf y)}{\xi_yF_{2}(\mathbf y)- \xi_{y^{\prime\prime}}F_{1}(\mathbf y)}&=&\frac{F_{2y^{\prime\prime}}(\mathbf y)-F_{3y^\prime}(\mathbf y)}{\xi_{y^{\prime\prime}}F_{2}(\mathbf y )-\xi_{y^\prime}F_{3}(\mathbf y)}
\end{eqnarray*} and they are functions in $\xi:=\beta(y)\gamma(y^\prime)\delta(y^{\prime\prime}).$
Moreover, the integrating factor is given by 
\[
\mu(\xi)=\exp\left\{\displaystyle\int\frac{F_{3t}(\mathbf y)-F_{0y^{\prime\prime}}(\mathbf y)}{\xi_{y^{\prime\prime}}F_{0}(\mathbf y)}d\xi\right\}.
\]
\end{corollary}
\section{Concluding Remarks}
In this paper, we investigated the existence of integrating factor of the following class of third order nonlinear differential equations: 
\begin{equation}\label{thirdorder12}
F_3(t,y,y^\prime,y^{\prime\prime})y^{\prime\prime\prime}+F_2(t,y,y^\prime,y^{\prime\prime})y^{\prime\prime}+F_1(t,y,y^\prime,y^{\prime\prime})y^\prime+F_0(t,y,y^\prime,y^{\prime\prime})=0.
\end{equation}
In fact, we presented some theoretical results related to the existence of certain forms of an integrating factor of \eqref{thirdorder12}. Also,  we presented some illustrative examples. In fact, these results not only useful for finding an integrating factor of \eqref{thirdorder12} analytically but also computationally. Particularly, we can check the validity of the conditions in our results by using the symbolic toolboxes in MATLAB and  MAPLE softwares.  Also, we can find the integrating factor using the formula in these results with the help of the symbolic toolboxes. Moreover, by using the same argument above, we can derive an integrating factor in terms of $\xi=\alpha(t)+\beta(y)+\gamma(y^\prime)+\delta(y^{\prime\prime})$.
 \section*{\textbf{Acknowledgment}} We would like to thank the editor and the referees for their valuable comments on this paper.   


\begin{thebibliography}{99}
\bibitem{AlAhmad} R. AlAhmad, M. Al-Jararha, H. Almefleh, \textit{Exactness of Second Order Ordinary Differential Equations and Integrating Factors,} JJMS \textbf{11}(3) (2015), pp. 155-167.

\bibitem{Aljararha} M. Al-Jararha, \textit{Exactness of Higher Order Nonlinear Ordinary Differential Equations}, Global Journal of Pure and Applied Mathematics, \textbf{11} (3) (2015), pp. 1467-1473.

\bibitem{existencepp} S. C. Anco and G. Bluman, \textit{Integrating factors and first integrals for ordinary differential equations},  Euro. Jnl of Applied Mathematics, \textbf{9} (1998), pp. 245-259. 

\bibitem{Ames}  W. F. Ames, \textit{Nonlinear ordinary differential equations in transport process}, Vol. \textbf{42}, Academic Press, New York, 1968.

\bibitem{Cheb}  E. S. Cheb-Terrab and  A. D. Roche , \textit{Integrating Factors for Second-order ODEs}, J. Symbolic Computation, \textbf{27} (1999),  pp. 501--519.

\bibitem{Davis} H. T. Davis, \textit{Introduction to nonlinear differential and integral equations}, Drover, New York, 1965.

\bibitem{kai}  K. Gehrs, \textit{Integrating Factors of Some Classes of Third Order ODEs}, Applied Mathematical Letters, \textbf{21} (2008),  pp. 748--753.

\bibitem{Jordan} D. W. Jordan and P. Smith, \textit{Nonlinear ordinary differential equations: An introduction for scientist and engineers}, 4th edition, Oxford University Press, 2007.

\bibitem{Bouquet} P. G. L. Leach and S. E. Bouquet, \textit{Symmetries and Integrating Factors}, J. Nonlinear Math. Phys., \textbf{9} (suppl. 2) (2002), pp. 73-91.

\bibitem{Lefschetz}   S. Lefschetz, \textit{Differential Equations: Geometric Theory}, Academic Press, New York, 1963.
\bibitem{cov} M. V. Makarets and V. Yu. Reshetnyak, \textit{Ordinary differential equations and calculus of variations}, World Scientific Publishing Co. Inc., NJ, 1995
\bibitem{Muriel} C. Muriel and J. L. Romero, \textit{First integrals, integrating factors and $\lambda$-symmetries of second-order differential equations},  J. Phys. A: Math. Theor., \textbf{42} (36) (2009):365207.





\bibitem{Peter}   P. J. Olver, \textit{Applications of Lie Group to Differential Equations}, Springer, USA, 1993.
\bibitem{Struble}  R. A. Struble, \textit{Nonlinear differential equations}, McGraw-Hill, New York, 1962.
\bibitem{von}  T. von Karman, \textit{The engineer grapples with nonlinear problems}, Bull. Amer. Math. Soc., \textbf{46} (1940),  pp. 615--683.
\end{thebibliography}
\end{document}